\documentclass[11pt]{article}
\usepackage{graphicx}
\usepackage{amssymb}

\def\C{\Gamma}
\def\Rbar{\overline{\bf R}}
\def\Rp{{\bf R}_{+}}
\def\R{{\bf R}}

\def\Chat{\hat{\bf C}}
\def\a{\alpha}
\def\b{\beta}
\def\d{\delta}
\def\e{\epsilon}

\def\o{\overline}

\def\li{{\rm li}}
\newtheorem{theorem}{Theorem}
\newtheorem{corollary}[theorem]{Corollary}
\newtheorem{lemma}[theorem]{Lemma}

\def\eproof{\rule{3mm}{3mm}\newline}

\begin{document}

\title{Distribution of intersection lengths  of a random geodesic with a geodesic lamination}
\author{Martin Bridgeman \\
 \footnotesize{Department of Mathematics,  Boston College,
Chestnut Hill, MA 02467}
\and
David Dumas\thanks{The second author was partially supported by
   an NSF postdoctoral research fellowship.} \\
\footnotesize{Department of Mathematics, Brown University, Providence, RI 02912}}
\date{November 10, 2006}

\maketitle

\begin{abstract}
We investigate the distribution of lengths obtained by intersecting a
random geodesic with a geodesic lamination. We give an explicit
formula for the distribution for the case of a maximal lamination and
show that the distribution is independent of the surface and
lamination.  We also show how the moments of the distribution are
related to the Riemann zeta function.
\end{abstract}

\section{Introduction}
Let $S$ be a closed hyperbolic surface and $\lambda$ a maximal
geodesic lamination on $S$. If we take an infinite geodesic $\a$ on
$S$, then $\lambda \cap \a$ decomposes $\a$ into geodesic arcs. We
consider the distribution of the lengths of these arcs. We show that
for almost every geodesic, the distribution is the same and we
explicitly calculate this distribution.

\begin{figure}
\begin{center}
\includegraphics[width=4in]{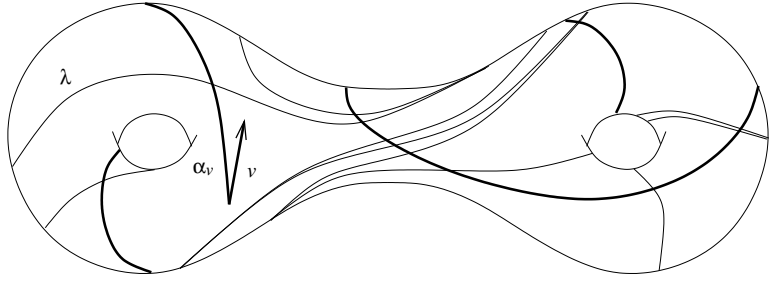}
\caption{Intersection of lamination $\lambda$ with ray $\a_{v}$ in direction $v$}
\label{setup}
\end{center}
\end{figure}

Let $v \in T_{1}(S)$ be a unit tangent vector on $S$ and $\a_{v}:\R
\rightarrow S$ the geodesic parameterized by arc length such that $v
= \a_{v}'(0)$. In the complement of $\lambda$, $\a_{v}$ is a countable
union of open intervals (see Figure \ref{setup}). We let
$$\a_{v}^{-1}(S - \lambda) = \cup_{j=1}^{\infty} I_{j}$$
where $I_{j}$ are open disjoint intervals. We denote the length $|I_{j}|$ of $I_{j}$ by $l_{j}$.

{\bf Definition:} The {\em length distribution} $D_{\lambda}(v)$ is a measure  on $(0, \infty)$ given by
\begin{eqnarray*}
D_{\lambda, t}(v) = \frac{1}{t}\sum_{I_{i} \cap [0,t] \neq \emptyset} \delta(l_{i}),\\
D_{\lambda}(v) = \overline{\lim_{t \rightarrow \infty}} D_{\lambda,t}(v)
\end{eqnarray*}
where $\d(x)$ is the Dirac measure at $x$ and the limit is taken with respect to the weak$^{*}$ topology.

The main result is that for almost all vectors $v$, the distribution $D_{\lambda}(v)$ is independent of $\lambda$ and $v$ and can be explicitly calculated.

\begin{theorem}
Let $\lambda$ be a maximal geodesic lamination. Then there exists a
measure $M$ such that for almost all  $v \in T_{1}(S)$ with respect to
the volume measure on $T_{1}(S)$ we have
$$D_{\lambda}(v) = \lim_{t \rightarrow \infty} D_{\lambda,t}(v) = M,$$
where 
$$dM = \frac{6xdx}{\pi^{2}\sinh^{2}{x}}.$$
\label{mainthm}
\end{theorem}

We define measure $P_{\lambda}(v)$ on $(0, \infty)$ similarly by
\begin{eqnarray*}
P_{\lambda, t}(v) = \frac{1}{t}\sum_{I_{i} \cap [0,t] \neq \emptyset} l_{i}\delta(l_{i}), \\
P_{\lambda}(v) = \overline{\lim_{t \rightarrow \infty}} P_{\lambda,t}(v).
\end{eqnarray*}

Note that if $\phi:(0, \infty) \rightarrow {\bf R}$ is a continuous function with compact support then $$P_{\lambda, t}(v)(\phi(x)) = D_{\lambda, t}(v)(x\,\phi(x)).$$
Then we obtain the following corollary:
\begin{corollary}
Let $\lambda$ be a maximal geodesic lamination. Then for almost all  $v \in T_{1}(S)$ with respect to the volume measure on $T_{1}(S)$, the measure $P_{\lambda}(v)$ is equal to the probability measure $P$ given by 
$$dP = \frac{6x^{2}dx}{\pi^{2}\sinh^{2}{x}}.$$
Furthermore the probability measure $P$  has expected value
$$E_{P}(x) = \frac{9}{\pi^{2}}\zeta(3),$$
and n$^{th}$ moment  given by
$$E_{P}(x^{n}) = \frac{3(n+2)!}{2^{n}\pi^{2}} \zeta(n+2),$$
where $\zeta$ is the Riemann zeta function. \end{corollary}

Before we give a proof of the corollary, let us describe a simple interpretation  of the probability measure $P$ in terms of the geometry of an ideal triangle. Consider picking a random tangent vector $p$ in an ideal triangle $T$ and let $g(p)$ be the geodesic through $p$. Then $P$ is the same probability distribution as the probability distribution of lengths of $T\cap g(p)$. In particular this corollary gives a nice description of $\zeta(3)$ in terms of the average length of a segment in  an ideal triangle.

{\bf Proof of corollary:}
Let $\phi:(0, \infty) \rightarrow {\bf R}$ be a continuous function with compact support.
Then by definition, $P_{\lambda, t}(\a_{v})(\phi(x)) = D_{\lambda, t}(v)(x\,\phi(x))$. Therefore as $x\,\phi(x)$ is continuous with compact support, it follows from the definition of the weak$^{*}$ topology that $P_{\lambda}(v)(\phi(x)) = D_{\lambda}(v)(x\,\phi(x))$. Therefore for almost all $v \in T_{1}(S)$,
$$P_{\lambda}(v)(\phi(x)) = \int_{0}^{\infty}x\,\phi(x) \:
 d(D_{\lambda}(v)) = \int_{0}^{\infty}\frac{6x^{2}\phi(x)\: dx}{\pi^{2}\sinh^{2}{x}}.$$
 It follows that for almost all $v \in T_{1}(S)$,
 $$d(P_{\lambda}(v)) = \frac{6x^{2} \:dx}{\pi^{2}\sinh^{2}{x}}.$$
We let $P$ be the measure on $(0,\infty)$ with distribution given above. To show that $P$ is a probability distribution and evaluate its moments, we calculate  the integral 
$$\int_{0}^{\infty}\frac{6x^{n}\: dx}{\pi^{2}\sinh^{2}{x}}$$
for $n \geq 2$. The antiderivative we describe below was found  using {\em Mathematica} (see 
\cite{Wolfram91}).

Let $\li_{n}(z)$ be the principal branch of the $n^{th}$ polylog function defined on the complex plane minus the real ray $[1,\infty)$.  The power series expansion about $0$ is given by
$$\li_{n}(x) = \sum_{k=1}^{\infty}\frac{x^{k}}{k^{n}}$$
and has radius of convergence $1$.

For $n=0$, $\li_{0}(x) = 1/(1-x)$, and for $n=1$,  $\li_{1}(x) = \log(1-x)$. Also for $n > 1$, we have convergence at $x=1$ and $\li_{n}(1) = \zeta(n)$. Taking the derivative of the power series we obtain the relation 
$$\li'_{n}(x) = \frac{\li_{n-1}(x)}{x}.$$ 
Let $n$ be a positive integer greater than 1. We define a function $F:(0,\infty) \rightarrow {\bf R}$ on the positive real axis  by
$$F(x) =- \sum_{k=0}^{n}\frac{n!}{k! 2^{n-k-1}}x^{k}\li_{n-k}(e^{-2x} ).$$
We now find the derivative of $F$. We have for $k \neq n$
$$\frac{d}{dx} x^{k}\li_{n-k}(e^{-2x}) =  kx^{k-1}\li_{n-k}(e^{-2x}) - 2x^{k}\li_{n-k-1}(e^{-2x}).$$
Also for $k= n$,  as $x^{n}\li_{0}(e^{-2x}) = x^{n}/(1-e^{-2x})$ we have 
$$\frac{d}{dx} x^{n}\li_{0}(e^{-2x}) =  nx^{n-1}\li_{0}(e^{-2x}) - \frac{2x^{n}e^{-2x}}{(1-e^{-2x})^{2}} = nx^{n-1}\li_{0}(e^{-2x}) - \frac{x^{n}}{2\sinh^{2}{x}}.$$
Putting all these derivatives together we have
$$F'(x) =  \frac{x^{n}}{\sinh^{2}{x}}.$$ 
Then 
$$ \int_{0}^{\infty}\frac{x^{n}\: dx}{\sinh^{2}{x}} = F(\infty) - F(0) = 0 - \left(-\frac{n!}{2^{n-1}}\li_{n}(1)\right) = \frac{n!}{2^{n-1}}\zeta(n).$$
Therefore as $\zeta(2) = \pi^{2}/6$
$$\int_{0}^{\infty}\frac{6x^{2}\: dx}{\pi^{2}\sinh^{2}{x}} = 1.$$
Therefore the n$^{th}$ moment of this probability distribution is
$$ \int_{0}^{\infty}\frac{6x^{n+2}\: dx}{\pi^{2}\sinh^{2}{x}} = \frac{3(n+2)!}{2^{n}\pi^{2}}\zeta(n+2). $$
In particular,  $P$ is a probability measure and has expected value $\frac{9}{\pi^{2}}\zeta(3)$.
\eproof

Given a geodesic current $m$, (a generalization of a closed geodesic described formally in the next section), we define a measure $D_{\lambda}(m)$ on $\R_{+}$ associated with the intersection of the geodesic current $m$ with the geodesic lamination $\lambda$. We prove the following properties of the function $m \rightarrow D_{\lambda}(m)$.

\begin{theorem}
Let $\mu$ be the Liouville geodesic current  for the closed hyperbolic surface $S$. Then 
\begin{enumerate}
\item $D_{\lambda}(\mu) = M$, where $M$ is as in Theorem \ref{mainthm}.
\item If $m_{i}$ are discrete geodesic currents such that $m_{i} \rightarrow \mu$ then
$$\lim_{i \rightarrow \infty} D_{\lambda}(m_{i}) = D_{\lambda}(\mu) = M.$$ 
\end{enumerate}
\label{main2}
\end{theorem}

Note that sequences $m_i$ as in part 2 above are abundant: Bonahon
shows in \cite{Bon88} that for almost any $v \in T_1(S)$, such a
sequence can be constructed from the geodesic ray $\alpha_v$
determined by $v$ (by closing up long segments).  The application of
Theorem \ref{main2} to this construction is discussed in Corollary
\ref{corbonahon} below.

\subsection*{Acknowledgments}
The authors thank the referee for comments and suggestions which
greatly shortened and simplified the paper.  The authors also thank Ed
Taylor and Curt McMullen for helpful discussions.

 \section{Geodesic laminations and geodesic currents}
 Let $S$ be a closed hyperbolic surface. A {\em geodesic lamination}
 is a closed subset of $S$ that is the union of a disjoint collection
 of geodesics on $S$ (see \cite{Thu79}). A geodesic lamination
 $\lambda$ is {\em maximal} if $S-\lambda$ is a disjoint collection of
 ideal triangles.
 
Let ${\bf H}^{2}$ be the hyperbolic plane and $G({\bf H}^{2})$ be  the set of oriented geodesics in  ${\bf H}^{2}$. We take the upper half plane model for ${\bf H}^{2}$ in the Riemann sphere $\Chat$. Then the boundary of ${\bf H}^{2}$ in $\Chat$ is $\Rbar = {\bf R} \cup\{\infty\}$ and is called the circle at infinity of ${\bf H}^{2}$. Identifying an oriented geodesic with its endpoints, we have that $G({\bf H}^{2}) = \Rbar \times \Rbar - \Lambda$ where $\Lambda$ is the diagonal in the product space. We give $G({\bf H}^{2})$ the subspace topology in $\Rbar \times \Rbar$.
 
In general, if $X$ is a topological space, we let $C_{0}(X)$  be the space of continuous functions on $X$ with compact support. If $m$ is a measure on $X$ we obtain  a linear function $m:C_{0}(X) \rightarrow {\bf R}$ by defining 

$$m(\phi) = \int_{X} \phi \: dm.$$

 The set of non-negative measures on $X$ is denoted $M(X)$. The weak$^{*}$ topology on $M(X)$  has basis at $m \in M(X)$ given by the sets $U(m, \phi, \e) = \left\{ n \in M(X) | \ | n(\phi) - m(\phi) | < \e\right\}$ where $\phi \in C_{0}(X), \e > 0$. Thus $m_{i} \rightarrow m$ if $m_{i}(\phi) \rightarrow m(\phi)$ for all $\phi \in C_{0}(X)$. If $A \subseteq X$ and $m \in M(X)$ we define $m|_{A}$ to be the restriction of the measure $m$ to $A$.

Let $S={\bf H}^{2}/\C$ be a closed hyperbolic surface, and $\alpha$  a closed geodesic in $S$ of length $l(\alpha)$. Then the  preimage of $\alpha$ in $G({\bf H}^{2})$ is a $\C$-invariant discrete subset. If $\a$ is primitive, we obtain a $\C$ invariant measure $m(\alpha)$ on $G({\bf H}^{2})$, by taking the Dirac measure on this set. If $\a$ is not primitive and $\a = \b^{k}$, where $\b$ is a primitive closed geodesic, we define $m(\a) = k\,m(\b)$.

This measure is an example of a {\em geodesic current}.

{\bf Definition:} A {\em geodesic current} on $S = {\bf H}^{2}/\C$ is a $\C$-invariant positive measure on $G({\bf H}^{2})$. The space of geodesic currents on $S$ is denoted ${\cal C}(S)$ and given the weak$^{*}$ topology.

If a geodesic current $m$ is a positive real multiple of $m(\alpha)$
for some closed geodesic $\alpha$ in $S$, then $m$ is called a {\em
discrete geodesic current}. If $m$ is a discrete geodesic current with
$m = \lambda.m(\a), \lambda \in {\bf R}_{+}$, we define the length
$l(m)$ of $m$ by $l(m) = \lambda l(\alpha)$ where $l(\alpha)$ is the
length of the closed geodesic $\alpha$.

Another geodesic current is the {\em Liouville} measure given by
$$\mu([a,b]\times[c,d]) = \left| \log \left| \frac{(a-c)(b-d)}{(a-d)(b-c)} \right| \right|.$$
By the invariance of the cross-ratio under M\"obius transformations,
$\mu$ is invariant under the group of isometries of ${\bf H}^{2}$. In
particular $\mu$ is invariant under the Fuchsian group $\C$ and is
therefore a geodesic current.  In fact, $\mu$ is the unique measure on $G({\bf H}^{2})$ (up to constant multiple) that is invariant under the full M\"obius group. Differentially $\mu$ is given by
$$d\mu_{(a,b)} = \frac{dxdy}{|a-b|^{2}}$$
where $dxdy$ is the standard area measure on ${\bf R} \times {\bf R}$. For a full description of geodesic currents and their properties see Bonahon \cite{Bon88}.

The following lemmas show that geodesic currents are a natural extension of closed geodesics.

\begin{lemma}[{Bonahon  \cite{Bon88}}]
 The set of discrete geodesic currents is dense in ${\cal C}(S)$. 
\end{lemma}

\begin{lemma}[{Bonahon  \cite{Bon88}}]
There exists a continuous function $l:{\cal C}(S) \rightarrow {\bf
  R}_{+}$ such that $l(m(\a)) = l(\alpha)$ for all closed geodesics
  $\alpha$, and $l(k\,m) = k\,l(m)$, for all $k \in {\bf
  R}_{+}, m \in {\cal C}(S)$.
\end{lemma}

\section{Ergodic theory}
Let $S$ be a closed hyperbolic surface and $\lambda$ a maximal lamination in $S$. Let $T_{1}(S)$ be the unit tangent bundle of $S$ and $\Omega$ the standard volume measure on $T_{1}(S)$. We normalize $\Omega$ to obtain the {\em unit volume measure} $V$ on $T_{1}(S)$. Then as $\Omega(T_{1}(S)) = 2\pi\:\mbox{Area}(S) = 4\pi^{2}|\chi(S)|$,
$$V = \frac{\Omega}{4\pi^{2}|\chi(S)|}.$$
 As before, if $v \in T_{1}(S)$, we let $\a_{v}:\R \rightarrow S$ be
 the geodesic in $S$ parameterized by arc length such that $\a'_{v}(0) = v$. Then $\a_{v}^{-1}(S-\lambda) = \cup_{i}I_{i}$ with $I_{i}$ having length $l_{i}$.  We define a length function $L_{\lambda}:T_{1}(S) \rightarrow [0,\infty]$ where $L_{\lambda}(v) = l_{j}$  if $\a_{v}(0) \in I_{j}$, and otherwise $L_{\lambda}(v) = 0$. Then $(L_{\lambda})_{*}V \in M(\Rp)$ is defined by $((L_{\lambda})_{*}V)(\phi) = V(\phi \circ L).$

\begin{theorem}
For almost all $v \in T_{1}(S)$ (wrt $\Omega$), $\lim_{t \rightarrow \infty} P_{\lambda, t}(v)$ exists and satisfies
$$\lim_{t \rightarrow \infty} P_{\lambda, t}(v) = P_{\lambda}(v) =(L_{\lambda})_{*} V.$$
\end{theorem}

{\bf Proof:}
Let $\phi \in C_{0}(\Rp)$, a continuous function on $\R_{+}$ with compact support. Then
$$P_{\lambda, t}(v)(\phi) = \frac{1}{t}\sum_{I_{i} \cap [0,t] \neq \emptyset} l_{i} \phi(l_{i}).$$ 

Let $L_{\lambda}:T_{1}(S) \rightarrow \R_{+}$ be the function
 described above. Then if $\a_{v}(t) \in I_{i}$, we have
 $L_{\lambda}(\a_{v}'(t)) = l_{i}$.  Therefore
  $$P_{\lambda, t}(v)(\phi) = \frac{1}{t}\sum_{I_{i} \cap [0,t] \neq
  \emptyset} l_{i} \phi(l_{i}) = \frac{1}{t} \sum_{I_{i} \cap [0,t]
  \neq \emptyset} \int_{I_{i}} \phi(L_{\lambda}(\a_{v}'(t))) \: dt =
  \frac{1}{t} \int_{J_{t}} \phi(L_{\lambda}(\a_{v}'(t))) \: dt$$
  where $J_{t}$ is the union of the intervals $I_{i}$ that intersect $[0,t]$. 
 We compare $P_{\lambda, t}(v)(\phi)$ to the time average, i.e.
$$Q_{\lambda, t}(v)(\phi) = \frac{1}{t}
 \int_{0}^{t}\phi(L_{\lambda}(\a_{v}'(s))) \: ds.$$
If $\a_{v}(0) \not\in \lambda$, then there exists an interval $I_{i_{0}}= (-a_0,b_0)$ with $0 \in I_{i_{0}}$.  If $\a_{v}(0)\in \lambda$, we set $a_0 = 0$.
Similarly if $\a_{v}(t) \not\in \lambda$, then there exists an interval $I_{i_{t}}= (t-a_{t},t + b_{t})$ with $t \in I_{i_{t}}$.  If $\a_{v}(t)\in \lambda$, we set $b_{t}= 0$.
Then
$$P_{\lambda, t}(v)(\phi) - Q_{\lambda, t}(v)(\phi) =
\frac{1}{t}\left( a_0 \phi(l_{i_{0}}) + b_{t}\phi(l_{i_{t}})\right).$$
As $\phi$ has compact support, let $\mbox{supp}(\phi) \subseteq
[x,y]$.  Then we have
\begin{equation}
|P_{\lambda, t}(v)(\phi) - Q_{\lambda, t}(v)(\phi)| \leq \frac{(a_0 + y) \:
 \| \phi \|_{\infty}}{t},
\label{PQcomp}
\end{equation}
where $\| \cdot \|_\infty$ is the sup norm.

The integral $Q_{\lambda,t}(v)(\phi)$ is the time average of the function $\phi \circ L_{\lambda}$  along the ray $t \mapsto \a_{v}'(t)$ in $T_{1}(S)$.
The function $\phi \circ L_{\lambda}$ is easily shown to be in
$L^{1}(V)$. As the geodesic flow on $T_{1}(S)$ is ergodic (see
\cite{Hop39}), by the Birkhoff ergodic theorem (see \cite{CFS82}) we have that for almost all $v$ with respect to the volume measure $V$, 
$$\lim_{t \rightarrow \infty}Q_{\lambda, t}(v)(\phi) = \int_{T_{1}(S) }\phi(L_{\lambda}(v)) \:dV = ((L_{\lambda})_{*}V)(\phi).$$
Therefore, for any $\phi \in C_{0}(\Rp)$ there exists a set of full measure $A_{\phi} \subset T_{1}(S)$ such that  $$\lim_{t \rightarrow \infty}Q_{\lambda,t}(v) (\phi) = ((L_{\lambda})_{*}V)(\phi) \mbox{ for all } v \in A_{\phi}.$$ 

By the Weierstrass approximation theorem (see \cite{Rud53}), we can
choose a countable basis $\{\phi_{i}\}$ for the sup norm topology on $C_{0}(\Rp)$. Then $A = \cap A_{\phi_{i}}$ is a set of full measure on which $\lim_{t \rightarrow \infty}Q_{\lambda,t}(v) (\phi_{i}) = ((L_{\lambda})_{*}V)(\phi_{i})$ for all $i$.

Now let $\phi \in C_{0}(\Rp)$ and $v \in A$. Then there exists a subsequence $\phi_{i_{j}}$ converging uniformly to $\phi$ on $\Rp$.

By linearity,  we have
$$|Q_{\lambda, t}(v)(\phi) - Q_{\lambda, t}(v)(\phi_{i_{j}})| \leq \|\phi_{i_{j}} - \phi\|_{\infty}\qquad
|((L_{\lambda})_{*}V)(\phi_{i_{j}}) - ((L_{\lambda})_{*}V)(\phi)| \leq \|\phi_{i_{j}} - \phi\|_{\infty}.$$

By uniform convergence of $\phi_{i_{j}}$, for any $\e > 0$, there exists an $m > 0$ such that $\|\phi_{i_{j}} - \phi \|_{\infty} < \e$ for all $j > m$. 
We choose a $j > m$, then
$$\begin{array}{ll}
|Q_{\lambda, t}(v)(\phi) - ((L_{\lambda})_{*}V)(\phi)| \leq & |Q_{\lambda, t}(v)(\phi) -Q_{\lambda, t}(v)(\phi_{i_{j}})| +  |Q_{\lambda, t}(v)(\phi_{i_{j}}) - ((L_{\lambda})_{*}V)(\phi_{i_{j}})| \\
&+ |((L_{\lambda})_{*}V)(\phi_{i_{j}}) - ((L_{\lambda})_{*}V)(\phi)|.
\end{array}
$$
Therefore 
$$|Q_{\lambda, t}(v)(\phi) - ((L_{\lambda})_{*}V)(\phi)| \leq   |Q_{\lambda, t}(v)(\phi_{i_{j}}) - ((L_{\lambda})_{*}V)(\phi_{i_{j}})| + 2\e.$$
As $v \in A$, $\lim_{t \rightarrow \infty}Q_{\lambda,t}(v)(\phi_{i_{j}}) = ((L_{\lambda})_{*}V)(\phi_{i_{j}})$. Thus there exists an $T > 0$ such that 
$$|Q_{\lambda,t}(v)(\phi_{i_{j}}) -
((L_{\lambda})_{*}V)(\phi_{i_{j}})| \leq \e \qquad \mbox{ for all } t \geq T.$$
Therefore
$$|Q_{\lambda,t}(v)(\phi) - (L_{\lambda})_{*}(V)(\phi)| \leq 3\e \qquad
\mbox{ for all } t \geq T$$
giving 
$$\lim_{t \rightarrow \infty}Q_{\lambda,t}(v)(\phi) = ((L_{\lambda})_{*}V)(\phi).$$
By (\ref{PQcomp}) there are constants $a_0, y$ such that
$$|P_{\lambda, t}(v)(\phi) - Q_{\lambda, t}(v)(\phi)| \leq \frac{(a_0 +
  y) \: \|\phi\|_{\infty}}{t}.$$
Therefore we have
$$\lim_{t \rightarrow \infty}P_{\lambda,t}(v)(\phi) = \lim_{t \rightarrow \infty}Q_{\lambda,t}(v)(\phi) = ((L_{\lambda})_{*}V)(\phi).$$
As $\phi$ is arbitrary,
$$\lim_{t \rightarrow \infty}P_{\lambda,t}(v) = (L_{\lambda})_{*}V \qquad \mbox{ for almost every } v \in T_{1}(S).$$
\eproof
 
Let $M({\bf R}_{+})$ be the space of measures on the open interval $(0,\infty)$ with the weak$^{*}$ topology. If $f:{\bf R}_{+} \rightarrow {\bf R}_{+}$ is a continuous function, then we obtain a map $\hat{f}:M({\bf R}_{+}) \rightarrow M({\bf R}_{+})$ given by 
$$\hat{f}(m)(\phi) = m(f\cdot\phi) = \int_{0}^{\infty} f(x) \phi(x) \: dm$$ for any $\phi \in C_{0}({\bf R}_{+})$ with compact support. On the level of infinitesimals, the map $\hat{f}$ is just multiplication by the function $f$ and thus has inverse  $\widehat{(1/f)}$. 

Taking $r(x) = 1/x$, we have by definition that 
$$\hat{r}(P_{\lambda, t}(v)) = D_{\lambda, t}(v).$$
Therefore we have the following corollary.
\begin{corollary}
For almost all $v \in T_{1}(S)$ (wrt $\Omega$), $\lim_{t \rightarrow \infty} D_{\lambda, t}(v)$ exists and satisfies
$$\lim_{t \rightarrow \infty} D_{\lambda, t}(v) = D_{\lambda}(v) = \hat{r}((L_{\lambda})_{*}V).$$
\label{dPisxdD}
\end{corollary}

From the above corollary, Theorem \ref{mainthm} follows from
calculating $(L_{\lambda})_{*}V$. Before we do the calculation, we
consider the simpler case of geodesics and geodesic currents
intersecting a single ideal triangle in ${\bf H}^2$.

\section{Intersection of a geodesic current with an ideal triangle}
We first consider the length of  intersection of a geodesic with a fixed ideal triangle. Let $T$ be an ideal triangle in ${\bf H}^{2}$. We define $L:G({\bf H}^{2}) \rightarrow [0,\infty]$ by 
 $L(g) = \mbox{Length}(g \cap T)$, with $L(g) = 0$ if $g \cap T = \emptyset$.  We observe that if $L(g)$ is non-zero and finite, then $L$ is continuous at $g$. Thus letting $O_{T} = L^{-1}((0, \infty))$ we have that $L$ is continuous when restricted to $O_{T}$.
    
If $m$ is a geodesic current,  we define a measure $L_{*}m$ on $(0,\infty)$ as follows. Let $\phi:(0,\infty) \rightarrow {\bf R}$ be a continuous map with compact support in $(0, \infty)$. We extend $\phi$ to a continuous map $\overline{\phi}:[0,\infty] \rightarrow {\bf R}$ by defining $\overline{\phi}(0) = \overline{\phi}(\infty) = 0$ and $\overline{\phi}(x) = \phi(x)$ otherwise.
 We then define  the measure $L_{*}m$ on $(0, \infty)$ by letting $(L_{*}m)(\phi) = m(\overline {\phi} \circ L)$. We will suppress the extension notation in what follows and identify $\phi$ with $\overline{\phi}$.
 
  We now describe the length function $L$ in terms of the endpoints of the geodesic. Let $T$ be the ideal triangle in the upper half space ${\bf H}^{2}$ with endpoints $\{0,1,\infty\}$ in $\Rbar$. We denote the intervals $I_{1} = (-\infty, 0)$, $I_{2} = (0,1)$ and $I_{3} = (1, \infty)$. We consider the  nine subsets $I_{i} \times I_{j} \subseteq G({\bf H}^{2})$. If $g \in I_{i} \times I_{i}$ then $g \cap T = \emptyset$ and therefore $L = 0$ on the subsets $I_{i} \times I_{i}$. If $g \not\in I_{i} \times I_{j}$ for any $i,j$, then $g$ has one endpoint in $\{0,1,\infty\}$ and therefore has either zero or infinite length. Therefore the set $O_{T} = L^{-1}((0,\infty))$ is given by
    $$O_{T} =   \bigcup^{3}_{\stackrel{i \neq j}{i,j}} I_{i} \times I_{j} .$$ 

We define $L_{ij}:I_{i}\times I_{j} \rightarrow (0, \infty)$ to be the restriction of $L$ to $I_{i}\times I_{j}$.

\begin{lemma}
The length function $L_{12}:I_{1}\times I_{2} \rightarrow (0,\infty)$ is given by the formula 
$$L_{12}(u,v) = \frac{1}{2}\ln\left(\frac{1-u}{1-v}\right).$$
\label{trianglelength}
\end{lemma}

\begin{figure}
\begin{center}
\includegraphics[width=4in]{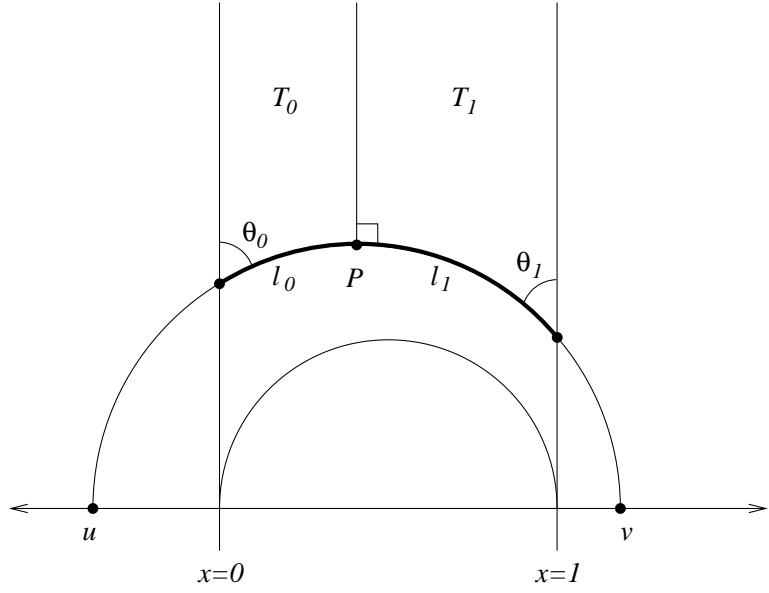}
\caption{The intersection of a geodesic $g$ with the
  standard ideal triangle $T$ in ${\bf H}^2$}\label{tintersect}
\end{center}
\end{figure}

{\bf Proof:} We will first find $L_{13}$. Let $g = (u,v) \in I_{1}
\times I_{3}$ be a geodesic. We drop a perpendicular $P$ from $\infty$
to $g$, with vertex $p$. We will first assume that $p \in T$. Then $P$
and $g$ decompose $T$ into two triangles an a quadrilateral. We label
the triangles $T_{0}$ and $T_{1}$ so that the vertical line $x=i$ is a
side of $T_i$.  Both triangles have a right angle and an ideal
vertex. We label the other angle of $T_{i}$ by $\theta_{i}$ and the
length of the only finite length side of $T_{i}$ by $l_{i}$ (see
Figure \ref{tintersect}). Then $L_{13}(u,v) = l_{0}+l_{1}$. By
hyperbolic trigonometry we have $\tanh{l_{i}} = \cos{\theta_{i}}$ (see
\cite{Thu79}).  We denote the euclidean center and radius of $g$ by
$c, r$ respectively. Then $u = c-r, v = c+r$. From euclidean
trigonometry we have
$$\tanh{l_{0}} = \cos{\theta_{0}} = \frac{c}{r} \qquad \tanh{l_{1}} = \cos{\theta_{1}} = \frac{1-c}{r}.$$
As $\tanh^{-1}{x} = \frac{1}{2}\ln(\frac{1+x}{1-x})$, we have
$$l_{0} = \frac{1}{2}\ln\left(\frac{1+\frac{c}{r}}{1-\frac{c}{r}}\right)= \frac{1}{2}\ln\left(-\frac{v}{u}\right), \qquad l_{1} = \frac{1}{2}\ln\left(\frac{1+\frac{1- c}{r}}{1-\frac{1-c}{r}}\right)= \frac{1}{2}\ln\left(\frac{-u+1}{v-1}\right).$$
Combining we obtain
$$L_{13}(u,v) = \frac{1}{2}\ln\left(\frac{v(u-1)}{u(v-1)}\right).$$
The case when $p \not\in T$ gives the same formula. 

The isometry $f(z) = 1/\overline{z}$ maps $I_{1}\times I_{3}$ to
$I_{1}\times I_{2}$. Therefore on $I_{1}\times I_{2}$ we have
$$L_{12}(u,v) = L_{13}(1/u, 1/v) =  \frac{1}{2}\ln\left(\frac{1-u}{1-v}\right).$$
\eproof

Recall that $\mu$ is the Liouville measure on $G({\bf H}^2)$.  We now use
Lemma \ref{trianglelength} to calculate its pushforward by $L$.

\begin{theorem}
 The measure $M_{T} = L_{*}\mu$ is given differentially by
 $$dM_{T} = \frac{6x\: dx}{\sinh^{2}{x}}.$$
\label{Liouvilleform}
 \end{theorem}

{\bf Proof:} As $\mu$ is invariant under M\"obius transformations, it
suffices to consider the standard ideal triangle $T$ with endpoints
$\{0,1,\infty\}$ in $\Rbar$.  As before we note that $L$ is supported
in $O_T = \cup_{i \neq j} I_i \times I_j$ where $I_{1} = (-\infty,
0)$, $I_{2} = (0,1)$ and $I_{3} = (1, \infty)$ as before note that

Let $\phi: (0,\infty) \rightarrow {\bf R}$ be a continuous function with compact support on $(0, \infty)$. Then  by definition of the pushforward, we have

$$(L_{*}\mu)(\phi) = \mu(\phi \circ L) =  \sum_{\stackrel{i \neq j}{i,j=1}}^{3}\left( \int_{I_{i}\times I_{j}} \phi \circ L \: d\mu \right) = \sum_{\stackrel{i \neq j}{i,j=1}}^{3} M_{ij}(\phi).$$

Thus $L_{*}\mu$ on $(0,\infty)$ is the sum of the six measures $M_{ij}, i \neq j$. 

Let $i \neq j$ and $k \neq l$, then there exists a M\"obius
transformation $p$ sending $T$ to itself and mapping $I_{i} \times
I_{j}$ homeomorphically to $I_{k} \times I_{l}$. By invariance of
$\mu$, we have $\mu = p_{*}\mu$.  The map $p$ is a hyperbolic isometry and sends $T$
to itself, thus $L \circ p = L$.  Therefore by change of variable in the integral $M_{ij}$ we have
$M_{kl} = M_{ij}$.  Thus
$$L_{*}\mu = 6\,M_{12}.$$
 
By Lemma \ref{trianglelength} we have that for $(u, v) \in I_{1}\times I_{2}$
$$L(u, v) =  \frac{1}{2}\ln \left( \frac{1- u}{1-v} \right). $$

 Letting $x = L(u,v)$, we solve for $u$ and $v$, to obtain $u = 1-
(1-v)e^{2x}$, $v = 1 -(1- u)e^{-2x}.$ The measure $\mu$ is given
differentially by
$$d\mu = \frac{dudv}{|u-v|^{2}}$$

Therefore 
$$M_{12}(\phi)  = \int_{I_{1}}\int_{I_{2}} \frac{\phi(L(u,v))\: dudv}{|u-v|^{2}}.$$
We make the change of variable $x = L(u,v)$, $v = v$. Then $u = f(v,x) = 1- (1-v)e^{2x}$ and $v = g(v,x) = v$. Thus the Jacobian is
$$J(v,x) =  \left| \begin{array}{cc} 
\frac{\partial f}{ \partial x} & \frac{\partial g}{ \partial x}\\
\frac{\partial f}{ \partial v}& \frac{\partial g}{ \partial v}
\end{array}\right| =
\left| \begin{array}{cc} 
-2(1-v)e^{2x} & 0\\
e^{2x} & 1
\end{array}\right|
= 2(1-v)e^{2x} 
.$$

Then we have
$$d\mu = \frac{J(v,x)\: dvdx}{|u-v|^{2}} = \frac{2(1-v)e^{2x}\: dvdx}{|1-(1-v)e^{2x}-v|^{2}} =  \frac{dvdx}{2(1-v)\sinh^{2}{x}}.$$

Let $F(u,v) = (v,x)$. Then we have that
$$M_{12}(\phi) =  \int_{F(I_{1}\times I_{2})} \frac{\phi(x)\: dvdx}{2(1-v)\sinh^{2}{x}}.$$

 We now need to find  $F(I_{1}\times I_{2})$. We fix a $v \in
 (0,1)$. Then from the formula for $L$ in terms of $u,v$ it follows that $L$ has range 
  $$ \frac{1}{2} \ln\left( \frac{1}{1-v}\right) <  L(u,v)   < \infty.$$
 Therefore
$$F(I_{1}\times I_{2}) = \left \{ (v,x) \:\left| \:0 < v < 1, \:
 \frac{1}{2} \ln\left( \frac{1}{1-v}\right) <  x   < \infty \right. \right \}$$
 rewriting we obtain
 $$F(I_{1}\times I_{2}) = \left \{ (v,x) \:| \:0 < v < 1-e^{-2x}, \: 0 <
 x   < \infty \right \}.$$
 
Therefore 
$$M_{12}(\phi) =  \int_{0}^{\infty} \left(\int_{0}^{1- e^{-2x}}\:
\frac{dv}{2(1-v)}\right)\frac{\phi(x)\: dx}{\sinh^{2}{x}}.$$
Substituting $t = 1-v$ we get

$$ \int_{0}^{1- e^{-2x}}\:\frac{dv}{2(1-v)} = \frac{1}{2} \int^{1}_{e^{-2x}}\:\frac{dt}{t} = \frac{1}{2}(\ln(1) - \ln(e^{-2x})) = x.$$
Therefore we have
$$M_{12}(\phi) =  \int_{0}^{\infty} \frac{x\phi(x)\:dx}{\sinh^{2}{x}},$$
giving
$$dM_{T}  = \frac{6x\:dx}{\sinh^{2}{x}}.$$
$\eproof$

Related to the measure $M_{T}$ is a probability measure $P_{T}$ on $(0,\infty)$ given by the pushforward of the unit volume measure on the unit tangent bundle $T_{1}(T)$ to an ideal triangle $T$. If $p \in T_{1}(T)$ then associated to $p$ we let $g(p)$ be the oriented geodesic with tangent vector $p$ . Then $L(g(p))$ is the length of the segment through $p$ in $T$. The measure $P_{T}$ is defined to be  the probability distribution of lengths for a randomly chosen tangent vector in $T$. 

If we let $\Omega$ be the standard volume measure on $T_{1}(T)$, then the probability measure $P_{T}$ is given by
$$P_{T} = \frac{(L\circ g)_{*}\Omega}{\Omega(T_{1}(T))}.$$

\begin{theorem}  The measure $P_{T}$ is given differentially by
$$dP_{T} =   \frac{6x^{2}dx}{\pi^{2}\sinh^{2}{x}}.$$
\label{Pcalc}
\end{theorem}

{\bf Proof:}
The map $g:T_{1}({\bf H}^{2}) \rightarrow G({\bf H}^{2})$ is a trivial fiber bundle. Thus we have $T_{1}({\bf H}^{2}) = G({\bf H}^{2}) \times {\bf R}$. In terms of this representation, the volume form $\Omega$ is given by $d\Omega = 2 d\mu \times dt$ where $\mu$ is the Liouville geodesic current and $dt$ is hyperbolic length along the fiber (see \cite{Nic89} for details).
Thus 
$$(L\circ g)_{*}\Omega(\phi) = \int_{T_{1}(T)} \phi(L(g(p)))\: d\Omega.$$
As $T_{1}({\bf H}^{2}) = G({\bf H}^{2}) \times {\bf R}$, we have $p$ represented by $(x,y,t)$ where $(x,y)$ are the endpoints of $g(p)$. Furthermore $d\mu = dxdy/|x-y|^{2}$. Also, for a given geodesic $(x,y)$, the parameter $t$ takes values from $c(x,y)$ to $c(x,y) + L(x,y)$, for some constant $c(x,y)$.Thus
$$(L\circ g)_{*}\Omega(\phi) = \sum_{\stackrel{i,j =1}{i \neq j}}^{3}
\int_{I_{i}\times I_{j}}\frac{2\: dxdy}{|x-y|^{2}}
\int_{c(x,y)}^{c(x,y) + L(x,y)}\phi(L(x,y)) \: dt.$$
Therefore integrating over $t$
$$(L\circ g)_{*}\Omega(\phi) = \sum_{\stackrel{i,j =1}{i \neq j}}^{3}
\int_{I_{i}\times I_{j}}\frac{2\phi(L(x,y))L(x,y)\:
  dxdy}{|x-y|^{2}}. $$
Using Theorem \ref{Liouvilleform} and performing the change of
variables $ l = L(x,y)$, $y  = y$, we obtain
$$(L\circ g)_{*}\Omega(\phi) = \sum_{\stackrel{i,j =1}{i \neq j}}^{3}
\int_{0}^{\infty} \frac{2l^{2}\phi(l)}{\sinh^{2}{l}}\: dl =
\int_{0}^{\infty} \frac{12l^{2}\phi(l)}{\sinh^{2}{l}}\: dl.$$
Therefore 
normalizing by $\Omega(T_{1}(T)) = 2\pi^{2}$, we obtain
$$dP_{T}(x) = \frac{6x^{2}}{\pi^{2}\sinh^{2}{x}}dx.$$
\eproof
An immediate corollary  relates the expected value of $P_{T}$ to the Riemann zeta function;
\begin{corollary}
The expected value $E(P_{T})$ of $P_{T}$ is  given by
$$E(P_{T}) = \int_{0}^{\infty} \frac{6x^{3}}{\pi^{2}\sinh^{2}{x}}\: dx =  \frac{9}{\pi^{2}}\zeta(3) = 1.09614\ldots$$
\end{corollary}

To put this expected value in perspective, we observe that it is quite
close to the diameter of the hyperbolic disk inscribed in $T$, which
is $\log(3) = 1.09861\ldots$.

\section{Calculating $(L_{\lambda})_{*}V$}
Let $S = {\bf H}^{2}/\C$ be a closed hyperbolic surface and let $\lambda$ be a maximal geodesic lamination on $S$. Then $S - \lambda = \cup_{i=1}^{k} T_{i}$ where $T_{i}$ are disjoint ideal triangles. Then as the area of $S$ is $2\pi|\chi(S)|$, we have $k = 2|\chi(S)|$.

We let $T'_{i}$ be an ideal triangle in ${\bf H}^{2}$ which is a lift of $T_{i}$. Then as before, we define
$L'_{i}(g) = \mbox{Length}(g\cap T'_{i})$. If $T''_{i}$ is another lift of $T_{i}$, then there exists a M\"obius transformation $p \in \C$ with $p(T'_{i}) = T''_{i}$. Let $L''_{i}(g) = \mbox{Length}(g \cap T''_{i})$. Then $L'_{i}(g) = \mbox{Length}(g \cap p^{-1}(T''_{i})) = \mbox{Length}(p(g) \cap T''_{i}) = L''_{i}(p(g))$. Therefore $L'_{i} = L_{i}''\circ p$. Thus  if $m \in {\cal C}(S)$ is a geodesic current, then $(L'_{i})_{*}m = (L_{i}'' \circ p)_{*}m = (L''_{i})_{*}(p_{*}m)$. As $m$ is invariant under $\C$, $p_{*}m = m$. Therefore $(L'_{i})_{*}m = (L''_{i})_{*}m$ and therefore can define $(L_{i})_{*}m = (L'_{i})_{*}m$.

{\bf Definition:} If $m$ is a geodesic current of length $l(m)$ then the {\em distribution of lengths of intersections of $m$ and $\lambda$}, denoted $D_{\lambda}(m)$, is the measure
$$D_{\lambda}(m) = \frac{1}{l(m)} \sum_{i=1}^{k} (L_{i})_{*}m .$$

From Theorem \ref{Liouvilleform} we have that on $(0,\infty)$
$$d(D_{\lambda}(\mu))  = \frac{12|\chi(S)| x\,dx}{l(\mu)\sinh^{2}{x}}.$$
 Then  as $l(\mu) = 2\pi^{2}|\chi(S)|$ (see Bonahon \cite{Bon88}), we have 
 $$d(D_{\lambda}(\mu)) = \frac{6x\:dx}{\pi^{2}\sinh^{2}{x}} = dM.$$ 
Thus $D_{\lambda}(\mu) = M$ and we call $M$
the  {\em distribution of lengths of intersections} for a maximal lamination. This proves part 1 of Theorem \ref{main2}.

\begin{theorem}
\label{pushforwardV}
For any maximal geodesic lamination $\lambda$, the pushforward of the
normalized volume measure $V$ on $T_1(S)$ by $L_\lambda$ is given by 
$$(L_{\lambda})_{*}V = P_{T}.$$
\end{theorem}

{\bf Proof:}
As $V$ is the unit volume measure obtained by normalizing $\Omega$,  and $\Omega(T_{1}(S)) = 4\pi^{2}|\chi(S)|$, we have

$$((L_{\lambda})_{*}V)(\phi) = \frac{1}{4\pi^{2}|\chi(S)|}\int_{T_{1}(S)} \phi(L_{\lambda}(v))\:d\Omega.$$
As $\lambda$ is a maximal lamination, we have $S-\lambda = \cup_{i=1}^{k} T_{i}$ where $T_{i}$ are disjoint ideal triangles. As $\lambda \subset S$ is measure zero, we have
$$((L_{\lambda})_{*}V)(\phi) = \frac{1}{4\pi^{2}|\chi(S)|}
\sum_{i=1}^{k} \int_{T_{1}(T_{i})} \phi(L_{\lambda}(v))\: d\Omega.$$

Let $q:T_{1}({\bf H}^{2}) \rightarrow T_{1}(S)$ be the standard
covering map. Let $T'_{i}$ be a lift of $T_{i}$ to ${\bf H}^{2}$. Then
we lift the integrals over $T_{i}$ using $q$ to integrals over
$T'_{i}$, also using $\Omega$ to denote the volume measure on
$T_{1}({\bf H}^{2})$, and obtain
$$ \int_{T_{1}(T_{i})} \phi(L_{\lambda}(v)) \: d\Omega =
\int_{T_{1}(T'_{i})} \phi(L_{\lambda}(q(v)) \: d\Omega.$$

Let $g:T_{1}({\bf H}^{2}) \rightarrow G({\bf H}^{2})$ be the fibre bundle such that $v$ is tangent to geodesic $g(v)$. As above, let $L_{i}:G({\bf H}^{2}) \rightarrow [0,\infty]$ where $L_{i}(g) = \mbox{Length}(g\cap T'_{i})$. Then for $v \in T_{1}(T'_{i})$,  $L_{\lambda}(q(v)) = L_{i}(g(v))$ and

$$  \int_{T_{1}(T'_{i})} \phi(L_{\lambda}(q(v)) \: d\Omega =
\int_{T_{1}(T'_{i})} \phi(L_{i}(g(v)) \: d\Omega.$$

By Theorem \ref{Pcalc},
$$ \int_{T_{1}(T'_{i})} \phi(L_{i}(g(v)) \: d\Omega =  2\pi^{2}  P_{T}(\phi).$$
Therefore as there are $2|\chi(S)|$ triangles $T_{i}$, we combine to obtain
$$((L_{\lambda})_{*}V)(\phi) =  \frac{1}{4\pi^{2}|\chi(S)|} \sum_{i=1}^{k} \left(2\pi^{2 }P_{T}(\phi)\right) = P_{T}(\phi).$$
Finally we have
$$d((L_{\lambda})_{*}V) = dP_{T} = \frac{6x^{2}\: dx}{\pi^{2}\sinh^{2}{x}}.$$
\eproof

Applying Corollary \ref{dPisxdD}, we obtain
\begin{corollary} For almost all $v \in T_{1}(S)$ (wrt $\Omega$), we
  have $D_{\lambda}(v) = M$ and thus
$$d(D_{\lambda}(v)) = \frac{6x\: dx}{\pi^{2}\sinh^{2}{x}}.$$
\end{corollary}

The above corollary completes the proof of Theorem 1.

\section{Continuity of $D_{\lambda}$ at $\mu$}
In this section we prove part 2 of of Theorem \ref{main2}. We recall the definition of $D_{\lambda}(m)$:
$$D_{\lambda}(m) = \frac{1}{l(m)} \sum_{i=1}^{k} (L_{i})_{*}m.$$

We first prove a lemma about weak$^{*}$ convergence.

\begin{lemma}
Let $X \subseteq {\bf R}^{n}$ be open and $K \subseteq X$ a compact domain with piecewise smooth boundary $\partial K$. If $m \in M(X)$ we define  $\overline{m} \in M(K)$ to be the restriction of $m$ to $K$. Let $m_{i} \rightarrow m$ in $M(X)$ and $\phi:K \rightarrow {\bf R}$  continuous  on $K$. Then
 $$\overline{\lim_{i\rightarrow \infty}} \left |\overline{m}_{i}(\phi)
 - \overline{m}(\phi) \right | \leq \int_{\partial K} |\phi| \: dm.$$ 
\label{technical1}
\end{lemma}

{\bf Proof:} 
Let  $\phi : K \rightarrow {\bf R}$ be a continuous  function. We let $\phi = \phi^{+} - \phi^{-}$ where  $\phi^{+}, \phi^{-}$ are both continuous non-negative functions on $K$. If the lemma is true for non-negative functions then we have
\begin{eqnarray*}
\overline{\lim_{i\rightarrow \infty}} \left |\overline{m}_{i}(\phi) -
\overline{m}(\phi) \right | & = & \overline{\lim_{i\rightarrow
    \infty}} \left |\overline{m}_{i}(\phi^{+}-\phi^{-}) -
\overline{m}(\phi^{+}-\phi^{-}) \right | \\
& \leq & \overline{\lim_{i\rightarrow \infty}}\left
|\overline{m}_{i}(\phi^{+}) - \overline{m}(\phi^{+}) \right | +
\overline{\lim_{i\rightarrow \infty}} \left
|\overline{m}_{i}(\phi^{-}) - \overline{m}(\phi^{-})\right | \\
& \leq & \int_{\partial K} \phi^{+} \: dm +  \int_{\partial K} \phi^{-} \: dm =  \int_{\partial K} |\phi| \: dm.
\end{eqnarray*}

We let $\phi : K \rightarrow {\bf R}$ be a continuous  non-negative function. Given an $\e > 0$, we let $N_{\e}(\partial K)$ be the $\e$ neighborhood of the boundary of $K$. For $\e$ small, we obtain  non-empty compact domains $K_{-\e}, K_{+\e}$ given by $K_{-\e} = K -  N_{\e}(\partial K)$, and $K_{+\e} = K \cup \overline{N_{\e}(\partial K)}$. Then $K_{-\e} \subset K \subset K_{+\e}$. We define continuous non-negative functions with compact support $\phi_{-\e}: X \rightarrow {\bf R}$, $ \phi_{+\e}: X \rightarrow {\bf R}$ which will approximate $\phi$ from above and below.

The function $\phi_{-\e}$ is chosen to be pointwise monotonically increasing in $\epsilon$ and satisfy
 $$\mbox{supp}(\phi_{-\e}) \subseteq K, \qquad \phi_{-\e}(x) = \phi(x) \mbox{ for } x \in K_{-\e}, \qquad \phi_{-\e}(x) \leq \phi(x) \mbox{ for } x \in K.$$

The function $\phi_{+\e}$ is similarly chosen to be pointwise monotonically decreasing in $\epsilon$ and satisfy
 $$\mbox{supp}(\phi_{+\e}) \subseteq K_{+\e},\qquad \phi_{+\e}(x) = \phi(x) \mbox{ for } x \in K, \qquad \|\phi_{+\e}\|_{X} \leq \|\phi\|_{K}.$$

Then as the measures are positive, we have
$$m_{i}(\phi_{-\e}) - m(\phi_{+\e}) \leq \o{m}_{i}(\phi) - \o{m}(\phi) \leq m_{i}(\phi_{+\e}) - m(\phi_{-\e}).$$
As $m_{i} \rightarrow m$ we take limits to obtain
$$m(\phi_{-\e}) - m(\phi_{+\e}) \leq \lim_{\o{i\rightarrow \infty}}(\o{m}_{i}(\phi) - \o{m}(\phi)) \leq  \o{\lim_{i\rightarrow \infty}}(\o{m}_{i}(\phi) - \o{m}(\phi)) \leq m(\phi_{+\e}) - m(\phi_{-\e}).$$
Thus letting  $\psi_{\e} = \phi_{+\e} - \phi_{-\e}$, we have
$$\overline{\lim_{i \rightarrow \infty}} \left |\o{m}_{i}(\phi) -
\o{m}(\phi) \right | \leq  \left |m(\psi_{\e}) \right |.$$

The function $\psi_{\e}$ is non-negative, monotonically decreasing in $\e$, and has $\mbox{supp}(\psi_{\e}) \subseteq \o{N_{\e}(\partial K)}$. Furthermore if  $x \not\in \partial K$,  then $\lim_{\e \rightarrow 0}\psi_{\e}(x) = 0$ and if $x \in \partial K$ then $\lim_{\e \rightarrow 0}\psi_{\e}(x) = \phi(x)$. Therefore by Lebesgue dominated convergence (see \cite{Rud53}) we have
$$\lim_{\e \rightarrow 0} m(\psi_{\e}) = \int_{\partial K} \phi \: dm. $$

Thus 
$$\overline{\lim_{i\rightarrow \infty}}|\o{m}_{i}(\phi) - \o{m}(\phi)| \leq \int_{\partial K} \phi \: dm. $$
\eproof

\begin{lemma}
Let $T$ be the lift of an ideal triangle embedded in $S$ and $\mu$ be the Liouville geodesic current of $S$. If $m_{i}$ is a sequence of discrete geodesic currents such that $m_{i} \rightarrow \mu$, then $L_{*}m_{i} \rightarrow L_{*}\mu$.
\label{lengthcont}
\end{lemma}

{\bf Proof:} We need to show that if $\phi:(0,\infty) \rightarrow {\bf
R}$ is continuous with compact support then $L_{*}m_{i}(\phi)
\rightarrow L_{*}\mu(\phi)$, or equivalently that $m_{i}(\phi \circ L)
\rightarrow \mu(\phi \circ L)$. This would follow directly from the
definition of the weak$^{*}$ topology if $\phi \circ L$ were
continuous and compactly supported, however this is not the case.  To
establish the convergence we will instead describe a continuous
function $\hat{L}$ such that $\phi \circ \hat{L}$ has compact support
and $m(\phi \circ L) = m(\phi \circ \hat{L})$. Then continuity follows
from the definition of the weak$^{*}$ topology.

Let $\mbox{supp}(\phi) = [a, b] \subseteq (0,\infty)$.  As before we
need only integrate $\phi \circ L$ over $O_{T} = L^{-1}( (0,\infty)
)$. We move $T$ by a M\"obius transformation to have endpoints
$\{0,1,\infty\}$ in the upper half space model of ${\bf H}^{2}$. Then
as in Lemma \ref{trianglelength}, the set $O_{T}$ is composed of 6
open sets $I_{i}\times I_{j}, i,j = 1,\ldots, 3, i\neq j$.

Let $A = L^{-1}([a, b])$ be the set of all geodesics in $O_{T}$ with length between $a$ and $b$. We let $A_{ij} = A \cap (I_{i}\times I_{j})$ and $L_{ij}:I_{i}\times I_{j} \rightarrow (0,\infty)$ be the restriction of $L$ to $I_{i}\times I_{j}$. Then for $m \in {\cal C}(S)$,

$$L_{*}m(\phi) = m(\phi \circ L) = \sum_{\stackrel{i\neq j}{i,i=1}}^{3} \left( \int_{I_{i}\times I_{j}} \phi \circ L_{ij}\:dm \right) = \sum_{\stackrel{i\neq j}{i,i=1}}^{3} m|_{I_{i}\times I_{j}}( \phi \circ L_{ij}).$$

By the formula for  $L_{12}$ we have if $(u,v) \in A_{12}$ then 
$$a \leq \frac{1}{2}\ln\left(\frac{1-u}{1-v}\right) \leq b.$$
Solving we get
$$(1-v)e^{2a} \leq 1-u \leq (1-v)e^{2b}.$$

As $v > 0$, we get $1-u \leq e^{2b}$ giving $1- e^{2b} \leq u$. Similarly as $u < 0$ we get $1  \leq (1-v)e^{2b}$ giving $v \leq 1- e^{-2b}$. Thus  we have $1- e^{2b} \leq u$ and $v \leq 1- e^{-2b}$.  Thus $A_{12} \subseteq [1- e^{2b},0) \times (0, 1- e^{-2b}]$.

We choose $\epsilon > 0$ so that if both $u,v$ are within $\e$ of $0$, then $L_{12}(u,v) < a$. Using the formula for $L_{12}$, we can choose 
 $\epsilon = (e^{2a}-1)/(e^{2a}+1)$.

Thus we consider the intervals 
$$J_{1} = [1-e^{2b}, -\e], J'_{1} = [-\e, 0] \mbox{ and } J_{2} = [\e, 1- e^{-2b}], J'_{2} = [0, \e].$$ Then the set $K = (J_{1} \times J_{2}) \cup (J_{1}\times J'_{2}) \cup (J'_{1} \times J_{2})$ is a compact subset of $G({\bf H}^{2})$. Also, we have that $\mbox{supp}(\phi \circ L_{12}) \subseteq K$. 

\begin{figure}
\begin{center}
\includegraphics[width=3in]{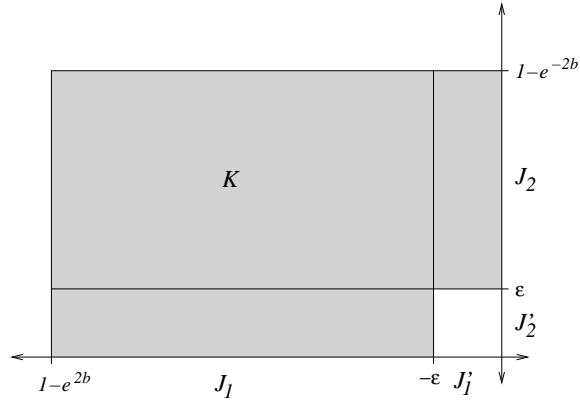}
\caption{The region $K$ as a subset of ${\bf R}^{2}$}
\end{center}
\end{figure}

The function $\phi \circ L_{12}$ is not continuous on $K$, being discontinuous at the points $D = (J_{1} \times\{0\}) \cup (\{0\} \times J_{2}) \subseteq \partial K$.  We define the function $\hat{L}_{12}: K \rightarrow {\bf R}$ by 
$$\hat{L}_{12}(u,v) = \frac{1}{2}\ln\left(\frac{1-u}{1-v}\right). $$
Note that this is the same formula as we obtained for $L_{12}$ in
Lemma \ref{trianglelength}, however, here we allow any $(u,v) \in K$ whereas
$L_{12}$ is supported on $I_1 \times I_2$.  The key difference is that
$K$ includes some geodesics (i.e. $D$) sharing one ideal endpoint with
$T$, and $L_{12}$ extends continuously to these points, giving
$\hat{L}_{12}$.

Thus $\hat{L}_{12}$ is continuous on $K$, with $\hat{L}_{12} = L$ on
$K-D$, and we have
$$\left| \left.m\right|_{K}(\phi \circ \hat{L}_{12}) - \left.m\right|_{I_{1} \times I_{2}}(\phi \circ L_{12}) \right|  = \left| \int_{D} \phi \circ \hat{L}_{12}  \: dm \right|  \leq m(D) \|\phi\|. $$

If $m = \mu$ then $m(D) = 0$ by the definition of $\mu$. As  $T$ is a lift of an embedded ideal triangle in $S$, none of the vertices of $T$ can be an endpoint of a closed geodesic on $S$. Therefore for $m$ discrete,   $m(D) = 0$. Thus for $m = \mu$ or $m_{i}$,
$$m|_{I_{1} \times I_{2}}(\phi \circ L_{12}) = m|_{K}(\phi \circ \hat{L}_{12}).$$
As $m_{i}\rightarrow \mu$, $K$ is a compact domain with piecewise smooth boundary, and, $\phi \circ \hat{L}_{12}$ is continuous, we have by Lemma \ref{technical1}:
$$\overline{\lim_{i\rightarrow \infty}}\left|\left.m_{i}\right|_{K}(\phi \circ \hat{L}_{12}) - \left.\mu\right|_{K}(\phi \circ \hat{L}_{12})\right| \leq \int_{\partial K} |\phi \circ \hat{L}_{12}| \: d\mu.$$
Therefore 
$$\overline{\lim_{i\rightarrow \infty}}\left| \left.m_{i}\right|_{I_{1} \times I_{2}}(\phi \circ L_{12}) - \left.\mu\right|_{I_{1} \times I_{2}}(\phi \circ L_{12})\right| \leq \int_{\partial K} |\phi \circ \hat{L}_{12}| \: d\mu.$$

We have $D \subseteq \partial K$ and by our choice of $K$ we have   $(\phi \circ \hat{L}_{12})(x) = 0$ for $x \in \partial K - D$. Therefore 

$$\overline{\lim_{i\rightarrow \infty}} \left| \left.m_{i}\right|_{I_{1} \times I_{2}}(\phi \circ L_{12}) - \left.\mu\right|_{I_{1} \times I_{2}}(\phi \circ L_{12}) \right| \leq \int_{D} |\phi \circ \hat{L}_{12}| \: d\mu \leq  \mu(D)\|\phi\|.$$
Therefore as $\mu(D) = 0$ we have
$$\lim_{i\rightarrow \infty} \left|\left.m_{i}\right|_{I_{1} \times I_{2}}(\phi \circ L_{12}) - \left.\mu\right|_{I_{1} \times I_{2}}(\phi \circ L_{12})\right| = 0.$$
Combining the contributions from each of the $I_{i}\times I_{j}$ we obtain 
$m_{i}(\phi \circ L) \rightarrow \mu(\phi \circ L).$
\eproof

We now use Lemma \ref{lengthcont} to prove the continuity property of $D_{\lambda}$ described in part 2 of Theorem~\ref{main2}.

\begin{corollary}
Let $\mu$ be the Liouville geodesic current of $S$. If $m_{i}$ is a sequence of discrete geodesic currents such that $m_{i} \rightarrow \mu$, then $D_{\lambda}(m_{i}) \rightarrow D_{\lambda}(\mu)$.
\label{contcor}
\end{corollary}

{\bf Proof:}
By definition
$$D_{\lambda}(m_{i})= \frac{1}{l(m_{i})} \sum_{j=1}^{k}
(L_{j})_{*}(m_{i}).$$ By Lemma \ref{lengthcont}, for each $j$ we have
$\lim_{i \rightarrow \infty} (L_{j})_{*}(m_{i}) =
(L_{j})_{*}(\mu)$. Also by continuity of the length function $l$, we
have $\lim_{i \rightarrow \infty} l(m_{i}) = l(\mu)$ and
$$\lim_{i \rightarrow \infty} D_{\lambda}(m_{i})= \lim_{i \rightarrow \infty}  \frac{1}{l(m_{i})} \sum_{j=1}^{k} (L_{j})_{*}(m_{i}) = \frac{1}{l(\mu)} \sum_{j=1}^{k} (L_{j})_{*}(\mu) = D_{\lambda}(\mu).$$
\eproof

We now apply Theorem \ref{main2} to Bonahon's construction of
sequences of discrete geodesic currents $m_i$ such that $m_i \to \mu$.

Recall that for any $v \in T_{1}(S)$ we defined the geodesic $\a_{v}$
by the map $\a_{v}:[0,\infty) \rightarrow S$ where $\a'(0) = v$ and
$\a$ is parameterized by arc length.  For any $L > 0$, we form a
closed path by joining the endpoints of $\a_{v}([0,L])$ by a shortest
arc (not necessarily unique).  If this curve is homotopic to a
geodesic we call this closed geodesic $\a^{L}_{v}$.

\begin{corollary}
\label{corbonahon}
For almost every $v \in T_{1}(S)$ with respect to the volume measure $\Omega$,  there exists a sequence $\{L_{n}\}$ monotonically increasing to infinity, such that  the geodesics $\b_{n} = \a_{v}^{L_{n}}$ satisfy
$$ \lim_{n \rightarrow \infty} D_{\lambda}(\b_{n})  = D_{\lambda}(\mu) = M.$$
\end{corollary}

{\bf Proof:} By Bonahon \cite{Bon88}, for almost every $v \in
T_{1}(S)$ with respect to the volume measure $\Omega$, the geodesics
$\a_{v}^{L}$ can be used to approximate the Liouville geodesic current
$\mu$. Specifically, for almost every $v \in T_{1}(S)$ there exists a
sequence $\{L_{n}\}$ monotonically increasing to infinity, such that
the geodesics $\b_{n} = \a_{v}^{L_{n}}$ exist and satisfy
$$\mu = (2\pi^{2} |\chi(S)|) \lim_{n \rightarrow \infty}  \frac{m(\b_{n})}{l(\b_{n})}.$$ 
Let $m_{i} = (2\pi^{2} |\chi(S)|)m(\b_{i})/l(\b_{i})$. By definition of $D_{\lambda}$ for geodesic currents, $D_{\lambda}(\b_{i}) = D_{\lambda}(m_{i})$. Therefore applying Corollary \ref{contcor} to the sequence $m_{i} \rightarrow \mu$, we have
$$ \lim_{i \rightarrow \infty} D_{\lambda}(\b_{i}) = \lim_{i \rightarrow \infty} D_{\lambda}(m_{i})  = D_{\lambda}(\mu) = M.$$
\eproof

\end{document}